\documentclass[12pt]{article}
\usepackage{amssymb}
\usepackage{amsthm}
\usepackage{graphicx}
\usepackage{graphics}
\usepackage{amsthm}

\input pictex.tex
\newcommand{\esp}{\hspace{0.05cm}}
\newcommand{\vsp}{\vspace{0.1cm}}

\newcommand{\clo}{\mathrm{S}^1}

\theoremstyle{definition}

\newtheorem{thm}{Theorem}[section]

\newtheorem{rem}[thm]{Remark}

\input epsf

\begin{document}

\date{}
\author{Andr\'es Navas}

\title{A finitely generated, locally indicable group with no 
faithful action by $C^1$ diffeomorphisms of the interval}
\maketitle

\vspace{-0.4cm}

\noindent{\bf Abstract.} According to Thurston's stability theorem, every group 
of $C^1$ diffeomorphisms of the closed interval is locally indicable ({\em i.e.}, 
every finitely generated subgroup factors through $\mathbb{Z}$). We show that, 
even for finitely generated groups, the converse of this statement is not true. 
More precisely, we show that the group $\mathbb{F}_2 \ltimes \mathbb{Z}^2$, 
although locally indicable, does not embed into $\mathrm{Diff}_+^1 (]0,1[)$. 
(Here $\mathbb{F}_2$ is any free subgroup of $\mathrm{SL}(2,\mathbb{Z})$, 
and its action on $\mathbb{Z}^2$ is the projective one.) Moreover, we show 
that for every non-solvable subgroup $G$ of $\mathrm{SL}(2,\mathbb{Z})$, the 
group $G \ltimes \mathbb{Z}^2$ does not embed into $\mathrm{Diff}^1_+(\clo)$.

\vspace{0.15cm}

\noindent{\bf MSC-class:} 20B27, 37C85, 37E05.

\vspace{0.75cm}

\noindent{\bf{\large{Introduction}}}


\vspace{0.5cm} Without any doubt, one of the most striking results about groups 
of diffeomorphisms is Thurston's stability theorem \cite{Th}. In the 
1-dimensional context, this theorem establishes that 
$\mathrm{Diff}^1_+([0,1[)$ is {\em locally 
indicable}, that is, each of its finitely generated subgroups factors through 
$\mathbb{Z}$. In the language of the theory of orderable groups, this is equivalent 
to saying that $\mathrm{Diff}^1_+([0,1[)$ is $C$-orderable (see for example 
\cite{order}). This is essentially the only known algebraic obstruction for 
embedding an abstract left-orderable group into $\mathrm{Diff}^1_+([0,1[)$.

A good discussion on {\em dynamical} obstructions for $C^1$ smoothability 
of continuous actions on the interval appears in D. Calegari's nice work 
\cite{calegari}. Most of them concern {\em resilient orbits}. Indeed, as was 
cleverly noticed by C. Bonatti, S. Crovisier, and A. Wilkinson, for groups of 
$C^1$ diffeomorphisms of the interval, there cannot be a central element without 
interior fixed points in the presence of resilient orbits 
\cite[Proposition 4.2.25]{book}. 
In the opposite direction, topologically transversal 
resilient orbits must appear when the topological entropy of the action 
is positive \cite{hurder}, or when some sub-pseudogroup acts without invariant 
probability measure \cite{DKN}. A new obstruction which does not involve resilient 
orbits is also given in \cite{calegari}. Nevertheless, these four conditions do 
not seem to complete the list of all possible dynamical obstructions. For 
instance, none of them seems to apply to groups of piecewise affine homeomorphisms, 
though `in general' the corresponding actions should be non $C^1$ smoothable... 

Giving a pure {\em algebraic} equivalent condition for the existence of 
a group embedding into $\mathrm{Diff}^1_+([0,1[)$ also seems very hard 
(see \cite{FF,growth} for two interesting particular cases). In this work, 
we show that local indicability, although necessary, 
is not a sufficient condition, even for finitely generated groups. 
For this, we deal with a concrete example, namely the group 
$\mathbb{F}_2 \ltimes \mathbb{Z}^2$ (which is easily seen to be 
locally indicable), where $\mathbb{F}_2$ is any free subgroup of 
$\mathrm{SL}(2,\mathbb{Z})$ whose action on $\mathbb{Z}^2$ 
is the projective one. 

\vspace{0.43cm}

\noindent{\bf Theorem A.} {\em The (locally indicable) 
group $\mathbb{F}_2 \ltimes \mathbb{Z}^2$ does not embed 
into $\mathrm{Diff}_+^1 (]0,1[)$.} 

\vspace{0.43cm}

The interest in considering the group $\mathbb{F}_2 \ltimes \mathbb{Z}^2$ 
comes from at least two sources. The first concerns the theory of orderable 
groups. Indeed, although $C$-orderable, this group admits no ordering with 
the stronger property of {\em right-recurrence}. This is cleverly noticed 
(and proved) in \cite{Witam}, where D. Witte-Morris shows that every 
finitely generated left-orderable amenable group admits a right-recurrent 
ordering, and hence every left-orderable amenable group is locally 
indicable. The second source of interest relies on Kazhdan's property (T). 
Indeed, from \cite[Th\'eor\`eme A]{superrig} it follows that, 
if the pair $(G,H)$ has the relative property (T) and $H$ is non-trivial and normal 
in $G$ \esp \esp 
-- as is the case of $(\mathbb{F}_2 \ltimes \mathbb{Z}^2, \mathbb{Z}^2)$ 
when $\mathbb{F}_2$ has finite index in $\mathrm{SL}(2,\mathbb{Z})$~--, 
\esp \esp then $G$ does not embed into the group of 
$C^{1 + \alpha}$ diffeomorphisms of the (closed) interval 
provided that $\alpha > 1/2$. It is perhaps possible to use the 
$L^p$ extensions of the (relative) property (T) in \cite{bader} to conclude, by 
a similar method, that $\mathbb{F}_2 \ltimes \mathbb{Z}^2$ does not embed into 
$\mathrm{Diff}_+^{1+\alpha}([0,1])$ for any $\alpha > 0$. However, it does not 
seem plausible to deal with the $C^1$ case (even for the closed interval) 
using this kind of arguments. (Algebraic obstructions for passing from 
$C^1$ to $C^{1+\alpha}$ embeddings exist: see for example \cite{growth}.) 

Our proof of Theorem A is strongly influenced by an argument due to 
J. Cantwell and L.~Conlon (namely the proof of the second half of 
Theorem 2.1 in \cite{CC}). It relies on considerations about `growth' of 
orbits (perhaps the right invariant to be considered should be the {\em topological 
entropy} associated to all possible actions on the interval). With slight modifications, 
these techniques also apply to the case of the circle. To motivate the theorem below, 
notice that $\mathrm{SL}(2,\mathbb{Z}) \ltimes \mathbb{Z}^2$ embeds into 
$\mathrm{Homeo}_+(\clo)$ (see \S \ref{ejemplos}).

\vspace{0.43cm}

\noindent{\bf Theorem B.} {\em For any non-solvable subgroup $G$ of 
$\mathrm{SL}(2,\mathbb{Z})$, the group $G \ltimes \mathbb{Z}^2$ does 
not embed into $\mathrm{Diff}_+^1(\clo)$.}

\vspace{0.43cm}
 
This result provides a first obstruction for group embeddings into 
$\mathrm{Diff}_+^1(\clo)$ for subgroups of $\mathrm{Homeo}_+(\clo)$ 
which does not rely on Thurston's stability theorem. This solves a question 
raised by J. Franks in a different manner from those of \cite{cal,parwani}.

\vsp

Unfortunately, our approach does not seem to be appropriate to deal with 
many other interesting groups which do act faithfully on the interval, 
as for example surface groups or general {\em limit groups} in the 
spirit of \cite{gel} (these groups are bi-orderable, which is 
stronger than being locally indicable). Another 
interesting question is the possibility of extending Theorem A 
to the group of {\em germs} of diffeomorphisms, where Thurston's 
theorem still applies (compare \cite[Remark 2.13]{growth}). Finally, 
the investigation of similar phenomena related to the higher 
dimensional versions of Thurston's theorem also seems promising.


\section{Existence of actions by homeomorphisms}
\label{ejemplos}

\hspace{0.45cm} As is well-known \cite{cal-book,book}, there exist 
faithful group actions of $\mathrm{SL}(2,\mathbb{Z}) \ltimes \mathbb{Z}^2$ 
by (orientation preserving) circle homeomorphisms. Indeed, let us consider 
the canonical action of $\mathrm{SL}(2,\mathbb{R})$ by real-analytic circle 
diffeomorphisms, and let $p \!\in\! \clo$ be a point whose stabilizer under this 
action is trivial. Replace each point $f(p)$ of the orbit of $p$ by an interval 
$I_f$ (where $f \in \mathrm{SL}(2,\mathbb{Z})$) in such a way that the total sum of 
these intervals is finite. Doing this, we obtain a topological circle $\clo_p$ provided 
with a faithful $\mathrm{SL}(2,\mathbb{Z})$-action (we use affine transformations 
for extending the maps in $\mathrm{SL}(2,\mathbb{Z})$ to the intervals $I_f$). 

Let $I \!=\! I_{id}$ be the interval corresponding to the point $p$, and let 
$\{\varphi^t \! : t \in \mathbb{R}\}$ be a non-trivial topological flow on $I$. 
Choose any real numbers $t_1,t_2$ which are linearly independent over the rationals, 
and let $h_1 = \varphi^{t_1}$ and $h_2 = \varphi^{t_2}$. Extend $h_1,h_2$ to $\clo_p$ 
by letting 
$$h_1(x) = f^{-1} \big( h_1^a h_2^c (f(x)) \big), \qquad 
h_2(x) = f^{-1} \big( h_1^b h_2^d (f(x)) \big),$$ 
where $x \!\in\! I_{f^{-1}}$ and 
\begin{equation}
f = \left(
\begin{array}
{cc}
a & b  \\
c & d  \\
\end{array}
\right) \esp \in \esp \mathrm{SL}(2,\mathbb{R}).
\end{equation}
For $x$ in the complement of the union of the $I_f$'s, we simply 
set $h_1(x)=h_2(x)=x$. The reader will easily check that the group 
generated by $\langle h_1,h_2 \rangle \sim \mathbb{Z}^2$ and the 
copy of $\mathrm{SL}(2,\mathbb{Z})$ acting on $\clo_p$ is 
isomorphic to $\mathrm{SL}(2,\mathbb{Z}) \ltimes \mathbb{Z}^2$.

If $\mathbb{F}_2$ is a free subgroup of finite index in $\mathrm{SL}(2,\mathbb{Z})$, 
then $\mathbb{F}_2 \ltimes \mathbb{Z}^2$ is locally indicable. Thus it acts faithfully 
by homeomorphisms of the interval \cite{order}. Although no such action arises as the 
restriction of the action constructed above, a faithful action may be constructed 
by following a similar procedure. For this, fix two (orientation 
preserving) homeomorphisms $f_1,f_2$ of 
$[0,1]$ generating a free group admitting a free orbit. There are many ways 
to obtain these homeomorphisms. We may take for example a left-ordering on 
$\mathbb{F}_2$, and next consider its dynamical realization (see the 
comment after Example 2.6 in \cite{order}). Another way is to use the 
fact that the group generated by $x \mapsto x+1$ and $x \mapsto x^3$ is 
free \cite{free}. Denoting by $p \!\in ]0,1[$ a point whose stabilizer 
under the corresponding $\mathrm{F}_2$-action is trivial, and then 
proceeding as above, we obtain the desired faithful action of 
$\mathrm{F}_2 \ltimes \mathbb{Z}^2$ on the interval.

Let us point out that, although the actions constructed above are only by homeomorphisms, 
they are topologically conjugate to actions by Lipschitz homeomorphisms 
(see \cite[Th\'eor\`eme D]{DKN}). 


\section{Preparation arguments: topological rigidity}
\label{top}

\hspace{0.45cm} Consider a faithful action of 
$\mathbb{F}_2 \ltimes \mathbb{Z}^2$ by homeomorphisms of the interval $[0,1]$. 
Let $I$ be an open (non-empty) {\em irreducible component} for the action of 
$\mathbb{Z}^2$, that is, a maximal open interval which is fixed by $\mathbb{Z}^2$. 
Since $\mathbb{Z}^2$ is normal in $\mathbb{F}_2 \ltimes \mathbb{Z}^2$, for every 
$f \in \mathbb{F}_2$ the interval $f(I)$ is also an open irreducible component 
for the action of $\mathbb{Z}^2$. 

\vspace{0.2cm}

According to \cite[\S 2.2.5]{book}, the group $\mathbb{Z}^2$ preserves a Radon measure 
$\mu$ on $I$. Associated to this measure, there is a non-trivial {\em translation number 
homomorphism} $\tau_{\mu} \!: \mathbb{Z}^2 \rightarrow \mathbb{R}$ defined by  
$\tau_{\mu} (g) = \mu ([x,g(x)[)$ for any $x \in I$. One has $\tau_{\mu} (g) > 0$ 
if and only if $g (x) \!>\! x$ for all $x \!\in\! I$. Moreover, if $\mu'$ is 
another invariant Radon measure, then $\tau_{\mu}$ and $\tau_{\mu'}$ coincide 
up to multiplication by a positive real number. 
We identify $h_1 \sim (1,0)$ and $h_2 \sim (0,1)$, 
and let $r = \tau_{\mu}((1,0))$ and $s = \tau_{\mu} ((0,1))$. 

\vspace{0.4cm}

\noindent{\bf Claim 1.} If $(r,s)$ is not an eigenvector of $f^T$, where 
$f \!\in\! \mathbb{F}_2$, then the interval $f(I)$ is disjoint from $I$.

\vspace{0.1cm}

\noindent{\em Proof.} Notice that $\tau_{\mu} (f(1,0)) = \tau_{f^*(\mu)} (1,0)$ and 
$\tau_{\mu} (f(0,1)) = \tau_{f^*(\mu)} (0,1)$. If $f$ fixes $I$, then $f^*(\mu)$ is 
another Radon measure on $I$ invariant by $\mathbb{Z}^2$. By the discussion above, 
there exists $\lambda > 0$ so that $\tau_{f^*(\mu)} = \lambda \tau_{\mu}$. 
This yields
$$\lambda r = \lambda \tau_{\mu}((1,0)) = \tau_{f^*(\mu)} ((1,0)) 
= \tau_{\mu} (f(1,0)) = \tau_{\mu} ((a,c)) = a r + c s.$$
Similarly, \esp $\lambda s = b r + d s.$ \esp 
This shows that $(r,s)$ is an eigenvector of $f^T$ with eigenvalue $\lambda$.

\vspace{0.5cm}

Now let $f_0$ be a hyperbolic element in $\mathbb{F}_2$ so that we have:

\vsp 

\noindent $(i)$ $\esp \esp$ $(r,s)$ is not an eigenvector of $f_0^T$, 

\vsp 
 
\noindent $(ii)$ $\esp \esp$ $(r,s)$ is not orthogonal to an eigenvector of $f_0^{-1}$, 

\vsp

\noindent $(iii)$ $\esp \esp$ neither $(1,0)$ nor $(0,1)$ are eigenvectors of $f$.

\vsp

\noindent By Claim 1, $f_0 (I)$ is disjoint from $I$. Thus, changing $f_0$ 
by its inverse if necessary, we may suppose that $f_0 (I)$ is to the left of 
$I$. Moreover, changing $f_0$ by $f_0^k$ 
for $k \!>\! 0$ sufficiently large, we may suppose 
that the expanding eigenvalue $\lambda$ of $f_0^{-1}$ is greater than 2. For 
a certain vector $(\alpha,\beta)$ in the expanding direction of $f_0^{-1}$ 
we have 
$$\lim_{n \rightarrow \infty} \big[ f_0^{-n} (1,0) - 
\lambda^n (\alpha,\beta) \big] = 0, \qquad 
\lim_{n \rightarrow \infty} \big[ f_0^{-n} (0,1) - 
\lambda^n (\alpha,\beta) \big] = 0.$$ 
The first of these equalities easily yields 
$$\lim_{n \rightarrow \infty} 
\big[ \tau_{\mu} (f_0^{-n} h_1 f_0^{n}) - \lambda^n (\alpha r + \beta s) \big] = 0.$$ 
Since $(r,s)$ is not orthogonal to any eigenvector of $f_0^{-1}$, the 
value of \esp \esp $t = \alpha r + \beta s$ is nonzero. \esp \esp 
Assume that $t$ is positive (the case where it is negative is similar). 
Replacing $h_1$ and $h_2$ respectively by $h_1^k$ and $h_2^k$ for $k > 0$ 
very large, we can ensure that $t > 0$ is sufficiently large so that 
we have: 

\noindent -- $\lambda t > 1$,

\noindent -- there exists an open interval $J \subset I$ with $0 < \mu (J) < t$, 

\noindent -- for all $i \in \mathbb{N}$ one has 
\begin{equation}
\label{secua}
i \leq 
t \left[ \lambda^i - \frac{\lambda^i - 1}{\lambda - 1} \right].
\end{equation}

\noindent Moreover, replacing $f_0$ by $f_0^k$ for $k > 0$ large 
enough, we may suppose that, for {\em all} $n \in \mathbb{N}$,  
\begin{equation}
\left| \tau_{\mu} (f_0^{-n} h_1 f_0^{n}) - \lambda^n t \right| 
\leq 1.
\label{casi-asimpt}
\end{equation}

Let $a$ (resp. $b$) be the fixed point of $f_0$ to the left (resp. to the right) 
of $I$. Since $f_0$ normalizes $\mathbb{Z}^2$, these points are also fixed by 
$\mathbb{Z}^2$. In \S \ref{half}, we will show that the dynamics of the subgroup 
$H$ of $\mathbb{F}_2 \ltimes \mathbb{Z}^2$ generated by $f_0$ and $h_1$ is not 
$C^1$-smoothable on $[0,1[$ by showing that, actually, it is not $C^1$-smoothable 
on $[a,b[$. The case of the open interval $]0,1[$ needs a supplementary 
argument and will be treated in \S \ref{open}.


\section{Cantwell-Conlon's argument: smooth rigidity}
\label{smooth}

\subsection{The case of the half-closed interval}
\label{half}

\hspace{0.45cm} In the statement of Cantwell-Conlon's theorem, there is an 
additional hypothesis of tangency to the identity at the endpoints. Nevertheless, 
such a hypothesis is not necessary, as the argument below shows.

\vspace{0.5cm}

\noindent{\bf Claim 2.} If the action of $\mathbb{F}_2 \ltimes \mathbb{Z}^2$ is 
by $C^1$ diffeomorphisms of $[0,1[$, then the restriction of $H$ to $[a,b[$ is 
topologically conjugate to a group of $C^1$ diffeomorphisms which are tangent 
to the identity at $a$.

\noindent{\em Proof.} This follows as a direct application of the M\"uller-Tsuboi's 
conjugacy trick: it suffices to conjugate by a $C^{\infty}$ diffeomorphism of $[a,b[$ 
whose germ at $a$ is that of $x \mapsto e^{-1/x^2}$ at the origin (see \cite{mu,Ts} 
for the details). 

\vspace{0.5cm}

In what follows, we will consider the dynamics of $f_0$ and $h_1$ after 
the preceding conjugacy, so they are tangent to the identity at $a$. 

\begin{rem} Since $h_1$ has a sequence of fixed points converging to $a$, its derivative 
at this point must equal 1 even for the original action; nevertheless, this was not 
necessarily the case for the original diffeomorphism $f_0$.
\end{rem}

\vspace{0.14cm}

\noindent{\bf Claim 3.} For each $k > 0$, the intervals of the form 
$$(f_{0}^{-k} h_1 f_0^k)^{\varepsilon_k} \cdots 
(f^{-2}_0 h_1 f_0^2)^{\varepsilon_2} (f^{-1}_0 h_1 f_0)^{\varepsilon_1}(J),$$
where $\varepsilon_i \in \{0,1\}$, are two-by-two disjoint.

\noindent{\em Proof.} Let 
$$W = (f_{0}^{-k} h_1 f_0^k)^{\varepsilon_k} \cdots 
(f^{-2}_0 h_1 f_0^2) (f^{-1}_0 h_1 f_0)^{\varepsilon_1}, \quad W' = 
(f_{0}^{-k} h_1 f_0^k)^{\varepsilon_k'} \cdots (f^{-2}_0 h_1 f_0^2)^{\varepsilon_2'} 
(f^{-1}_0 h_1 f_0)^{\varepsilon_1'}$$  
be such that $W \neq W'$, where all $\varepsilon_i,\varepsilon_i'$ 
belong to $\{0,1\}$. Let $i$ be the largest index for which $\varepsilon_i 
\neq \varepsilon_i'$, say $\varepsilon_i = 1$ and $\varepsilon_i' = 0$, 
and let   
$$W_* = (f_0^{-i} h_1 f_0^i) (f_0^{-(i-1)} h_1 f_0^{i-1})^{\varepsilon_{i-1}} 
\cdots (f_0 h_1 f_0)^{\varepsilon_1}, \quad 
W_*' = (f_0^{-(i-1)} h_1 f_0^{i-1})^{\varepsilon_{i-1}'} 
\cdots (f_0^{-1} h_1 f_0)^{\varepsilon_1'}.$$ 
Notice that each of the maps 
$(f_0^{-j} h_1 f_0^j)^{\varepsilon_j}$ either fixes all the points 
in $I$ (when $\varepsilon_j = 0$) or moves all of them to the right (when 
$\varepsilon_j = 1$). In particular, $W_*$ moves the left endpoint $u$ of 
$J \!= ]u,v[$ to a point $u_*$ which coincides with or is to the right of 
$f_0^{-i} h_1 f_0^i (u)$. By (\ref{casi-asimpt}), this implies that
\begin{equation}
\mu ([u,u'[) \geq \mu ([u,f_0^{-i} h_1 f_0^i (u)[) = 
\tau_{\mu} (f_0^{-i} h_1 f_0^i) \geq \lambda^i t - 1.
\label{prima}
\end{equation}
On the other hand, $W_*'$ moves $v$ to a point $v_*'$ which coincides with or 
is to the left of 
$$(f_0^{-(i-1)} h_1 f_0^{i-1}) \cdots (f_0^{-1} h_1 f_0)(v).$$
Since $\mu(J) < t$ and 
\begin{eqnarray*}
\mu \big( [v,(f_0^{-(i-1)} h_1 f_0^{i-1}) \cdots (f_0^{-1} h_1 f_0)(v)[ \big) 
&=& \tau_{\mu} \big( (f_0^{-(i-1)} h_1 f_0^{i-1}) \cdots (f_0^{-1} h_1 f_0) \big)\\ 
&=& \sum_{j=1}^{i-1}\tau_{\mu} (f_0^{-j} h_1 f_0^j)\\ 
&\leq& \sum_{j=1}^{i-1} \Big( \lambda^j t + 1 \Big)\\ 
&=& \left[ \frac{\lambda^i - 1}{\lambda - 1} - 1 \right] t + (i-1),
\end{eqnarray*}
inequalities (\ref{secua}) and (\ref{prima}) show that $v_*'$ is to the left 
of $u_*$. This implies that $W_*(J)$ and $W_*'(J)$ do not intersect, and 
hence $W(J) \cap W'(J) = \emptyset$.

\vspace{0.5cm}

To conclude the proof of the fact that the action of 
$\mathbb{F}_2 \ltimes \mathbb{Z}^2$ is not by $C^1$ 
diffeomorphisms of $[0,1[$, fix $N \!\in\! \mathbb{N}$ so 
that, for all $x \in [a,b[$ to the left of $f^N (I)$,
$$f_0'(x) \geq \sqrt[3]{3/4}, \quad 
\quad h_1'(x) \geq \sqrt[3]{3/4}.$$ 
Fix also a positive lower bound $A < 1$ for 
the derivative of $f_0$ and $h_1$ to the left of $I$. 
By opening brackets in the next expression, one easily 
checks that the length of each interval of the form  
$$(f_{0}^{-k} h_1 f_0^{k})^{\varepsilon_{k}} \cdots 
(f^{-2}_0 h_1 f_0^2)^{\varepsilon_2} (f^{-1}_0 h_1 f_0)^{\varepsilon_1}(J)$$
is at least
$$A^{3N} \Big( \sqrt[3]{\frac{3}{4}} \Big)^{3(k - N)} |J|.$$
Since there are $2^{k}$ of these intervals, 
for some constant $C > 0$, this yields
$$\left| [a,b] \right| \geq 
C \esp \Big( \frac{3}{2} \Big)^k |J|.$$
However, this is clearly impossible for a large $k$, thus completing the proof.

\vsp\vsp

We close this Section by noticing that similar arguments to those above apply 
to actions by $C^1$ diffeomorphisms of the interval $]0,1]$ instead of $[0,1[$. 


\subsection{The case of the open interval}
\label{open}

\hspace{0.45cm} To prove Theorem A in the general case of the {\em open} interval, 
we would like to apply the arguments of the preceding Section. For this, we need 
to ensure that either $a$ or $b$ actually belongs to $]0,1[$. Indeed, if not, 
we are not allowed to use the procedure of the Claim 1. 

Thus, we need to find a hyperbolic element $f \in \mathbb{F}_2$ such that:

\vsp

\noindent $(i)$ \esp \esp $(r,s)$ is not an eigenvector of $f^T$. 

\vsp

\noindent $(ii)$ \esp \esp $(r,s)$ is not 
orthogonal to any eigenvector of $f^{-1}$,

\vsp 

\noindent $(iii)$ $\esp \esp$ neither $(1,0)$ nor $(0,1)$ are eigenvectors of $f$, 

\vsp

\noindent $(iv)$ \esp \esp $f$ has fixed points inside $]0,1[$.

\vsp

\noindent For this, we begin by noticing that our free group 
$\mathbb{F}_2 \subset \mathrm{SL}(2,\mathbb{Z})$ must contain a free 
subgroup $F$ on two generators whose non-trivial elements are hyperbolic and 
satisfy properties $(i)$, $(ii)$, and $(iii)$ above. Indeed, this can be 
easily shown using a ping-pong type argument on $\mathbb{RP}^1$. Now  
$F$ must contain non-trivial elements having fixed points in $]0,1[$;  
if not, the action of $F$ on $]0,1[$ would be free, which is 
in contradiction with H\"older's theorem \cite{ghys,book}. Therefore, 
any element $f \in F$ having fixed points in $]0,1[$ satisfies $(i)$, 
$(ii)$, $(iii)$, and $(iv)$, and this concludes the proof of Theorem A.


\subsection{The case of the circle}
\label{general}

\hspace{0.45cm} Let $G$ be a non-solvable subgroup of 
$\mathrm{SL}(2,\mathbb{Z})$. To show Theorem B, we would again 
like to apply similar arguments to those of \S \ref{half}. However, 
there are certain technical issues that need a careful treatment. 

First of all, notice that, {\em a priori}, an irreducible component $I$ for the action of 
$\mathbb{Z}^2$ is not necessarily an interval: it could coincide with the whole circle. We 
claim, however, that this cannot happen. Indeed, let $(r',s')$ be the point in $\mathbb{T}^2$ 
whose coordinates are the rotation numbers of $(1,0)$ and $(0,1)$, respectively. Recall that 
the rotation number function is invariant under conjugacy. Moreover, its restriction to 
$\mathbb{Z}^2$ is a group homomorphism into $\mathbb{T}^2$ (see for example 
\cite[\S 6.6]{ghys} or \cite[\S 2.2.2]{book}). Since $G$ normalizes 
$\mathbb{Z}^2$, for all $f \!=\! \left(
\begin{array}
{cc}
a & b  \\
c & d  \\
\end{array}
\right) \! \in G$ we have (modulo $\mathbb{Z}$)
$$r' = \rho(1,0) = \rho(f(1,0)) = \rho(a,c) = a r' + c s'$$
and 
$$s' = b r' + d s'.$$ 
This means that $(r',s')$ is a fixed point for the action of $f^T$ on 
$\mathbb{T}^2$. But since $G$ is non-solvable, this cannot hold for 
every $f \in G$, thus showing that $I$ does not coincide with $\clo$. 

Now let $\mu$ be a $\mathbb{Z}^2$-invariant Radon measure on $I$. Let 
$\tau_{\mu} \! : \mathbb{Z}^2 \rightarrow \mathbb{R}$ be the corresponding 
translation number homomorphism, and let $r = \tau_{\mu}((1,0))$ and 
$s = \tau_{\mu}((0,1))$. Analogously to the case of \S \ref{open}, 
we need to find a hyperbolic element $f \in G$ so that the 
following conditions are fulfilled:

\vsp

\noindent $(i)$ \esp \esp $(r,s)$ is not an eigenvector of $f^T$,  

\vsp

\noindent $(ii)$ \esp \esp $(r,s)$ is not 
orthogonal to any eigenvector of $f^{-1}$,

\vsp

\noindent $(iii)$ $\esp \esp$ neither $(1,0)$ nor $(0,1)$ are eigenvectors of $f$, 

\vsp

\noindent $(iv)$ \esp \esp $f$ has fixed points on the circle.

\vsp

\noindent To obtain the desired element, we need to consider two cases separately. 

\vsp

If $G$ does not preserve any probability measure on $\clo$, then Margulis' alternative 
\cite{margulis} and its proof provide us with a free subgroup (in two generators) 
of $G$ whose elements have fixed points. Clearly, some (actually, `most') 
of these elements are hyperbolic and satisfy the conditions $(i)$, $(ii)$, 
and $(iii)$ above.

If $G$ preserves a probability measure on $\clo$, then the rotation number function 
$\rho \!: G \rightarrow \mathbb{T}^1$ is a group homomorphism \cite{ghys,book}. 
Therefore, the rotation number of all of the elements in $[G,G]$ is zero, and 
hence these elements must have fixed points. Since $G$ is non-solvable, $[G,G]$ 
contains free subgroups, which allows arguing as in the preceding case.


\vspace{0.5cm}

\noindent{\bf Acknowledgments.} The author is indebted to C. Bonatti, D. Calegari, 
L. Conlon, and T.~Gelander, for motivating and useful discussions on the subject. 

This work was funded by the PBCT-Conicyt Research Network on Low Dimensional Dynamics 
and the Math-Amsud Project DySET.


\begin{small}

\vspace{0.1cm}

\noindent Andr\'es Navas\\

\noindent Dpto. de Matem\'aticas, Fac. de Ciencia, Univ. de Santiago de Chile  

\noindent Alameda 3363, Santiago, Chile\\

\noindent E-mail address: andres.navas@usach.cl\\

\end{small}

\end{document}